\documentclass[11pt]{article}
\usepackage{amsfonts,amsthm,amsmath,amscd,fancyheadings,floatflt}
\divide\oddsidemargin by 2
\theoremstyle{plain}
\newtheorem{theorem}{\bf Theorem}
\newtheorem{corollary}{\bf Corollary}
\newtheorem{lemma}{\bf Lemma}
\newtheorem{proposition}{\bf  Proposition}
\newtheorem*{thmA}{\bf Theorem A}
\newtheorem*{thmB}{\bf Theorem  B}
\newtheorem*{coroA}{\bf Corollary A}

\theoremstyle{definition}
\newtheorem*{definition}{\bf Definition}

        {\end{list}}
 
\newsavebox{\savepar}

\newcommand{\C}{{\mathbb C}}
\newcommand{\N}{{\mathbb N}}

\newcommand{\Z}{{\mathbb Z}}

\newcommand{\chat}{{\widehat \C}}

\begin{document}

\title{Ruelle operator and transcendental entire maps}

\author{P. Dom\'{\i}nguez, P. Makienko and G. Sienra} 
\date{}

\maketitle

\begin{abstract}

If $f$ is a transcendental entire function with only algebraic singularities 
we calculate the Ruelle operator of $f$. Moreover, we prove both  (i) if $f$
 has a summable critical point, then $f$ is not structurally stable under 
certain topological conditions and (ii) 
if all critical points of $f$  belonging to Julia set  are summable, then   
there exists  no invariant lines fields in the Julia set.

\renewcommand{\thefootnote}{}
\addtocounter{footnote}{-2}           
\end{abstract}      


\section{Introduction}

If $f$ is a transcendental  entire map  we denote by  
$f^n$, $n \in \N$ , the  n-th iterate of $f$ and write the Fatou set as 
$F(f) = \{z \in \C$; there is some open set $U$  containing $z$ in which 
$\{f^n \}$ is a normal family $\}$. The 
complement of $F(f)$ is called  the Julia set $J(f)$. We  say that 
$f$ belongs to the class $S_q$ if the set of  singularities of $f^{-1}$ 
contains at most $q$ points.  Two entire maps $g$ and $h$ 
are topologically equivalent if there exist homeomorphisms 
$\phi, \psi: \C \to \C$ such that  $\psi \circ g = h \circ \phi$.

If we denote  by $M_f$, $f \in  S_q$ the set of all 
entire maps  topologically equivalent to $f$  we can define on $M_f$ as in 
\cite{eremenko} a structure of $(q + 2)$ dimensional complex manifold.

Fatou's conjecture states that the only structurally stable maps on $M_f$
are the hyperbolic ones. This conjecture is false in the case when there is 
an invariant line field in the Julia set of $f$, our result is a partial 
answer to this conjecture  for transcendental entire maps  with only 
finite number of  algebraic singularities.

P. Makienko \cite{Mak1, Mak2} and G.M. Levin \cite{lev} have  studied the
 Ruelle operator and the invariant line fields for rational maps,  the idea 
of this work is to study  an application
of the proposed approach  given in \cite{Mak2}  for  transcendental entire 
functions in class $S_q$,  where the singularities of $f^{-1}$ are only 
algebraic.  
\par{\bf Assumptions on maps. }
From now on we will assume that 
\begin{enumerate}
\item $f$ is transcendental entire and that the 
singularities of $f^{-1}$ are algebraic and finite and all critical points are simple (that is $ f''(c) \neq 0$). 

\item It follows from a very well known result of complex variables that there exist 
a decomposition
$$
\frac{1}{f'(z)}= \sum_{i=1}^{\infty} \left( \frac{b_i}{z-c_i}-p_{i}(z) \right)
+h(z),
$$
where $p_{i}$ are polynomials, $h(z)$ is an entire function, $\{c_{i}\}$ are the critical points of $f$, and $ b_{i} = \frac{1}{f''(c_i)}$ are constants depending on $f$. Now we assume that the series 
$$
\sum_{i=1}^{\infty}  \frac{b_i}{c_i^3}
$$
is absolutely convergent.
\end{enumerate}
\par Note that elements of generic subfamily  of the family $ P_1(z) + P_2(sin(P_3(z))) $ satisfy to assumptions above, here $ P_i(z) $ are polynomials.

Let $F_{n,m}$ the space of forms of the kind $\phi(z)Dz^mD\overline{z}^n$. Consider two formal  actions of $f$ on $F_{n,m}$,  say $f^*_{n,m}$ and 
${f_*}_{n,m}$,  on a function $\phi$ at the point $z$ by the formulas
\begin{eqnarray}
f^*_{n,m}(\phi) & = & \sum\phi(\xi_i)(\xi_i')^n(\overline{\xi_i'})^m = 
\sum_{y \in f^{-1}(z)}\frac{\phi(y)}{(f'(y))^n(\overline{f'(y)})^m}, \nonumber
\end{eqnarray}
and 
\begin{eqnarray}
{f_*}_{n,m}(\phi) & = & \phi(f)\cdot (f')^n\cdot\overline{(f')^m}, \nonumber
\end{eqnarray}

where $n, m \in \Z$ and $\xi_i, i = 1, ..., d$ are the branches of the 
inverse map $f^{-1}$. As in \cite{Mak1} we define  

\begin{enumerate}
\item The operator $f^* = f^*_{2,0}$ as {the \rm Ruelle operator} of 
      $f$.

\item  The operator $\vert f^*\vert = f^*_{1, 1}$ as  {\rm the modulus of the 
Ruelle operator}.

\item  The operator  $B_f = {f_*}_{-1,1}$ as the {\rm Beltrami operator}       of $f$.
\end{enumerate}

Let $c_i $ be the critical points of $f$ and  $Pc(f) = \overline{\cup_i\cup_{n \geq 0}f^n(f(c_i))}$ be the postcritical set. 

\begin{lemma}
Let $ Y \subset \chat $ be completely invariant measurable subset respect
 to $f$. Then
\begin{enumerate}
\item $ f^* : L_1(Y) \rightarrow L_1{Y} $ is linear endomorphism ``onto'' with
$\Vert f^*\Vert_{L_1(Y)} \leq 1$;
\item  Beltrami operator $ B_f: L_\infty(Y)\rightarrow L_\infty(Y) $ is dual operator to $f^*$;
\item if $ Y \subset \left \{ \chat \backslash \overline{\cup_if^{-i}(Pc(f))}\right \} $ is an open subset and let $ A(Y) \subset L_1(Y) $ be subset of holomorphic functions, then $ f^*(A(Y)) \subset  A(Y)$;
\item fixed points on the modulus of Beltrami operator define a non- negative
 absolutely  continuous invariant measure in $\C$.
\end{enumerate}
\end{lemma}

\noindent Observe that  all items above  follow from definitions.

\begin{definition} The space of quasi-conformal
deformations of a given map $f$, denoted by $ qc(f)$, is defined as.
\begin{eqnarray}
qc(f) = \bigl\{ g \in M_f: \text{ there is a quasiconformal
automorphism $ h_g $ of the Riemann} \nonumber \\
\text{ sphere } \chat  \text{ such that }
 g = h_g \circ f \circ {h^{-1}_g}\bigr\}{\big/} A_{\rm ff} (\C),  \nonumber
\end{eqnarray}
where $A_{ff} (\C)$ is the affine  group.
\end{definition}

\begin{definition} 
For $f \in S_q$ structurally stable the space of all grand orbits of 
$f$ on $\chat \setminus \overline{ \{ \cup_i f^{i}(P_c(f)) \} }$ forms  a 
disconnected Riemann surface, say $S(f)$,  of finite quasi-conformal type, 
see for details \cite{MS}.  
\end{definition}

\begin{definition} 
A  point $a \in f$ is called "{\rm summable}" 
if and only if either 
\begin{enumerate}
\item the set $X_a (f) = \overline{ \{ \cup_n f^n(f(a)) \} }$ is bounded and the series
$$\sum_{i = 0}\frac{1}{(f^i)\prime(f(a))}  
$$ 
is absolutely convergent or 
\item the set $X_a(f)$ is unbounded and the series
$$\sum_{i = 0}\frac{1}{(f^i)\prime(f(a))}  \text{ and }  \sum_{i = 0}\frac{\vert f^n(f(a)) \vert \vert \ln \vert f^n(f(c)) \vert \vert}{(f^i)\prime(f(a))}
$$ 
are absolutely convergent.

\end{enumerate}
\end{definition}

\begin{definition} Let $X$ be the space of transcendental entire  maps 
$f \in S_q$, fixing $0,1$,  with summable critical point $c \in J(f)$  
and either
\begin{enumerate}
\item $f^{-1}(f(c))$ is not in $X_c(f)$,  

\item $X_c(f)$ does not separate the plane, 
 
\item $m(X_c) = 0$, where $m$ is the Lebesgue  measure,  

\item $c \in \partial D \subset J(f)$, where $D$ is a component of 
$F(f)$.

Note that (4) includes the maps with completely invariant domain.
\end{enumerate}
\end{definition}

The main results of this work are Theorems A and B for transcendental 
entire maps. In \cite{Mak2} the theorems  were proved  for rational 
maps. The big differences 
between them is that for transcendental entire  maps infinity is an 
essential singularity and there are not poles.

\begin{thmA} Let $f \in X$. If $f$ has a summable 
critical point, then $f$ is not structurally stable map.
\end{thmA}


\begin{definition} Denote by  
$W$ the space of transcendental entire maps in $S_q$ such that:
\begin{enumerate} 
\item There is no parabolic points for $f \in W$.
\item All critical point are simple (that is $ f''(c) \neq 0$) and the forward orbit of any critical point $c$ is infinite and does not intersect the forward orbit of any other critical point.
\item $f$ satisfies (1) to (5) in the above definition, for all critical 
points of $f$.
\end{enumerate} 
\end{definition}

Conditions (1) and (2) are required for simplicity of the proof 
but they are not relevant.

\begin{definition} We call a transcendental  entire map $f$ {\rm summable} 
if all critical points  belonging to the Julia set are summable. 
\end{definition}
 
\begin{thmB} If $f \in W$ is  summable, then there exists no invariant line 
fields on $ J(f).$
\end{thmB}

\noindent {\bf Remark:} A theorem of McMullen \cite{McM} states that for 
the full family  $f_{a,b} = a + b\sin z$ has $m(J(f)) > 0$, then the arguments 
of J. Rivera Letelier \cite{Let} make non sense in this case .

\subsection*{Acknowledgements}
The authors would like to thank CONACYT  and the seminar of Dynamical 
Systems. This work was partially supported by proyecto CONACyT 
$\#$ 27958E, $\#$ 526629E  and  UNAM grant PAPIIT $\#$ IN-101700.

\section{Bers map}

Let $ \phi \in L_\infty(\C) $ and let $B_f(\phi) = \phi(f)\frac{{\overline f'}}{f'} : L_\infty(\C) \rightarrow L_\infty(\C) $ be the Beltrami operator. Then the open unit ball $ B $ of the space $ Fix(B_{f}) \subset L_\infty(\C) $ of fixed points of $ B_f$ is called {\it the space of Beltrami differentials for $f$} and describe all quasi-conformal  deformations of $f$.
\par Let $ \mu \in Fix(B_{f}), $ then for any $ \lambda $ with $ \vert\lambda\vert < \frac{1}{\Vert\mu\Vert} $ the element $ \mu_\lambda = \lambda\mu \in B\subset Fix(B_{f}$. Let $ h_\lambda $ be quasi-conformal maps corresponding to Beltrami differentials $ \mu_\lambda $ with $ h_\lambda(0, 1, \infty) = (0, 1, \infty).$ Then the map

$$ \lambda \rightarrow f_\lambda = h_\lambda\circ f\circ h_\lambda^{-1} \in M_{f} 
$$

\noindent is a conformal map. If  $ f_\lambda(z) = f(z) + \lambda G_\mu(z) + ..., $ then differentiation respect to $ \lambda $ in the point $ \lambda = 0$ gives the following equation

$$ F_\mu(f(z)) - f'(z)F_\mu(z) = G_\mu(z),
$$

\noindent where $ F_\mu(z) = \frac{\partial f_\lambda(z)}{\partial\lambda}\vert_{\lambda = 0}.$
\vspace{.5cm}

\noindent {\bf Remark 1.} Due to quasiconformal  map theory (see for example 
\cite{Krush}) for any $ \mu \in L_\infty(\C) $ with $\Vert\mu\Vert_\infty < \epsilon $ and small $ \epsilon$, there exists the following formula for quasi-conformal $ f_\mu $ fixing $ 0, 1, \infty.$

$$ 
f_\mu(z) = z -\frac{z(z - 1)}{\pi}\iint_{\C}\frac{\mu(xi) d\xi \wedge d\bar{\xi} }{\xi(\xi - 1)(\xi - z)} + C(\epsilon, f)\Vert\mu\Vert_\infty^2,
$$

\noindent where $ \vert z\vert < f$ and $ C(\epsilon, f) $ is constant does not depending on $ \mu.$ Then

$$
F_\mu(z) = \frac{\partial f_\lambda}{\partial\lambda}_{\vert\lambda = 0} = - \frac{z(z - 1)}{\pi}\iint_{\C}\frac{\mu(xi) d\xi \wedge d\bar{\xi}}{\xi(\xi - 1)(\xi - z)}.
$$

Hence we can define the linear map $ \beta: Fix(B_{f})\rightarrow H^1(f) $ by the formula, where $H^1(f)$ is defined below

$$
 \beta(\mu) = F_\mu(f(z)) - f'(z)F_\mu(z).
$$ 
 
We  call $ \beta$  {\it the Bers map } as an  analogy with Kleinian group (see for example \cite{ Kra}).
\bigskip

Let $ A(S(f))$ be the space of quadratic holomorphic integrable differentials on disconnected surface $ S(f)$.  Let $HD(S(f))$ be the space of harmonic differentials on $S(f))$. In every chart every element $\alpha \in HD(S(f))$ has a form $\alpha = \frac{\overline{\phi} d z^2}{\rho^2\vert dz\vert^2}, $ where $\phi d z^2 \in A(S(f))$ and $\rho\vert d z\vert$ is the Poincare metric. Let $P :\overline{\chat} \setminus \{ \overline{\cup_i f^{-i}(Pc(f))} \} \rightarrow S(f)$ be the projection. Then the pull back $ P_*: HD(S(f))\rightarrow Fix(B_{f}) $ defines  a linear  injective map. 

 The space $ HD(f) = P_*(A(S(f))) $ is called {\it the space of harmonic differentials.} For any element $ \alpha \in HD(f) $ the support $ supp(\alpha) \in F(f)$. Then $ dim(HD(f)) = dim(A(S(f))).$

Let $ J_f = Fix(B_{f})_{\vert J(f)} $ be the space of invariant Beltrami differentials supported by Julia set.
\bigskip

Now define  $H(f) = \{ \varphi: \Delta \to C_f  \text{ such that } 
\varphi(\lambda ) = f + \lambda f_1 + \lambda^{2} f_2 + ..., \text{ for } \lambda \text{ very small} \}$, and $C_f = \{g \in M_f : g(0) = 0 \text{ and } g(1) = 1 \} \subset M_f$.
\bigskip

We can define  an  equivalence 
relation $\sim$ on $H(f)$  in the following way, $\varphi_1 \sim \varphi_2$ 
if and only if $\frac{\partial (\varphi_1 - \varphi_2 )}{\partial \lambda} \vert_{\lambda = 0} = 0$.

\begin{definition} $H^1(f) = H(f)/ \sim$.
\end{definition}

Observe that (i) $H^1(f)$ is linear complex space and (ii) there exists an  
injection $\Psi$ such that $\Psi: H^1(f) \to T_f(C_f) =$ Complex tangent space.

\begin{theorem} Let $ f $ be structurally stable transcendental entire map. 
Then $\beta: HD(f)\times J_f \rightarrow H^1(f) $ is an isomorphism.
\end{theorem}

In structurally unstable cases $ \beta $ restricted on $ HD(f)\times J_f $ is always injective.
\bigskip

{\it Proof.} The map $f$ is structurally stable hence 
$dim(qc(f)) = dim(HD(f)\times J_f) = dim(H^1(f)) = dim(M_f /A_{\rm ff}(\C)) = q$. If we show that $\beta$ is {\it onto},  then we are done. 

Let $f_1$ be any element of $ H^1(f)$. There exists  a function 
$\varphi(\lambda)$ such that  for $\lambda$ sufficiently small  
$\varphi(\lambda) \subset C_f$ (since $f$ is structurally stable). Then  
$\varphi(\lambda) =  f_\lambda$ is a  holomorphic family of 
transcendental entire maps, thus 
$f_\lambda = h_\lambda\circ f\circ h_\lambda^{-1}$,  where $h_{\lambda}$ is 
a holomorphic family of quasi-conformal maps. Hence
  
$$f_1(z) = V(f(z)) - f'V(z),$$

where 
$V = \frac{\partial h_\lambda}{\partial \lambda} |_{\lambda = 0}$. 
The family of the complex dilatations $ \mu_\lambda(z) = \frac{\overline{\partial}h_\lambda(z)}{\partial h_\lambda(z)} \in Fix(f) $ forms a meromorphic family of Beltrami differentials. If $ \mu_\lambda(z) = \lambda\mu_1(z) + \lambda^2\mu_2(z) + ..., $ where $ \mu_i(z) \in Fix(f). $  Then 
$$
\frac{\partial h_\lambda}{\partial\lambda}_{\vert\lambda = 0} = -\frac{z(z - 1)}{\pi}\iint\frac{\mu_1(\xi) d\xi d\bar{\xi}}{\xi(\xi - 1)(\xi - z)} = F_{\mu_1}(z)
$$
and hence $ F_{\mu_1}(f(z)) - f'(z)F_{\mu_1}(z) = f_1$.
 
\bigskip

If we let  $ \nu = {\mu_1}_{\mid F(f)}$, then we can state the 
following claim.
\bigskip

\noindent {\bf Claim.} {\it There exists an element}  $ \alpha \in HD(S(f)) $
{\it such that} $ \beta(\alpha) = \beta(\nu). $
\bigskip

{\it Proof of the claim.} We will use here quasi-conformal theory  (see for 
example the  books of I. Kra \cite{Kra} and S.L. Krushkal \cite{Krush} and 
the papers of C. McMullen and D. Sullivan \cite{McM},  \cite{MS}). Let 
$\omega $ be the Beltrami differential on $ S(f) $ generated by $ \nu $ 
(that is $ P_*(\omega) = \nu$). Let $ <\psi, \phi> $ be the Petersen scalar 
product on $ S(f), $ where $ \phi, \psi \in A(S(f)) $ and

$$
  <\psi, \phi>  = \iint_{S(f)} \rho^{-2}\overline{\psi}\phi,
$$

\noindent where $ \rho $ is hyperbolic metric on disconnected surface $ S(f). $
Then by (for example) Lemmas 8.1 and 8.2 of chapter III in \cite{Kra} this scalar product defines a Hilbert space structure on $ A(S(f))$.
Then there exists an element $ \alpha^\prime \in HD(S(f)) $ such that equality 
$$ \iint_{S(f)}\omega\phi = \iint_{S(f)}\alpha^\prime\phi
$$
holds for all $ \phi \in A(S(f)).$

 Now let $ A(O) $ be space of all holomorphic integrable functions over $ O $, where $ O = \{\overline{\C}\backslash \overline{\cup_if^{-i}(Pc(f))}\} \subset F(f). $ Then the push forward operator $ P^*: A(O)\rightarrow A(S(f)) $ is dual to the pull back operator $ P_*.$ Hence element $ P_*(\alpha^\prime) $ satisfies the next condition
$$\iint_{O}\nu g = \iint_{O}P_*(\alpha^\prime) g, 
$$
for any $ g \in A(O). $

 All above means that $ \iint P_*(\alpha)\gamma_a(z) = \iint\nu\gamma_a(z) $ for all $ \gamma_a(z) = \frac{a(a -1)}{z(z - 1)(z - a)}, a \in J(f). $ Hence the  transcendental entire maps  $ \beta(P_*(\alpha))(a) = \beta(\nu)(a) $ on $ J(f) $ and we have the desired result with $ \alpha = P_*(\alpha^\prime). $ Thus the claim and the theorem are proved.


\section{Calculation of the Ruelle operator} 

Let us recall that from above there exist 
a decomposition
\begin{equation}
\frac{1}{f'(z)}= \sum_{i=1}^{\infty} \left( \frac{b_i}{z-c_i}-p_{i}(z) \right)
+h(z),
\end{equation}
where $p_{i}$ are polynomials, $h(z)$ is an entire function, $\{c_{i}\}$ are the critical points of $f$, and $ b_{i} =\frac{1}{f''(c_i} $ and the series $ \sum\frac{b_i}{c_i^3} $ is absolutely convergent.

In order to use Bers' density theorem and the infinitesimal formula of quasi-conformal maps, see Remark 1.  We  will work with linear combinations of the 
following functions.

$${\gamma_{a}}(z)=\frac{a(a-1)}{z(z-1)(z-a)} \in L_1(\C),$$ 
where $a \in \C \setminus \{0,1 \}$.

\begin{proposition}
Let ${\gamma_{a}}(z)$ as above. If $f$ is  any transcendental entire map with 
simple critical points, then 

\begin{displaymath}
f^{*}({\gamma}_{a}(z))= \frac{{\gamma}_{f(a)}(z)}{f'(a)}+\sum_{i=1}b_{i}{\gamma}_{a}(c_i){\gamma}_{f(c_{i})}(z).
\end{displaymath}
The coefficients $b_i$ and $c_i$ comes from (1).
\end{proposition}

{\it Proof.}  Let ${\varphi} \in C^{\infty}(S)$, where 
$S = \C \setminus \{0,1 \}$,  with 
compact support, denoted by $supp({\varphi})$ and 
${\varphi}(0)={\varphi}(1)=0$. Now consider the following:

\begin{displaymath}
\int_{\C}\varphi_{\bar z} f^{*}{\gamma_{a}}(z)dz \wedge d{\bar z}= \int_{\C}\frac {(\varphi_{\bar z} \circ f) {\bar f'}}{f'(z)}{\gamma_{a}}(z)dz \wedge d{\bar z}=
\int_{\C}\frac{(\varphi \circ f)_{\bar z}}{f'(z)}{\gamma_a}(z)dz \wedge d{\bar z}
\end{displaymath}
the first equality is by the duality with the Beltrami operator, see Lemma 1 in Section 1.

Let us denote by ${\psi}= {\varphi}(f)$, so $supp({\psi})= f^{-1}supp({\varphi})$ and is the union $\bigcup K_i$ of compact sets if there is not asymptotic 
values on it. Hence applying the decomposition in (1) of $1/f'(z)$ we have

\begin{gather}
\int_{\C}\frac{(\varphi \circ f)_{\bar z}}{f'(z)}{\gamma_a}(z)dz \wedge d{\bar z}=
\sum \int_{\C}{\psi}_{\bar z}\frac{a(a-1)b_i}{(z-c_{i})z(z-1)(z-a)}dz \wedge d{\bar z}-\nonumber\\
- \sum \int_{\C}{\psi}_{\bar z}p_{i}(z){\gamma_{a}}(z)dz \wedge d{\bar z}
+\int_{\C}{\psi}_{\bar z}h(z){\gamma_{a}}(z)dz \wedge d{\bar z}= \nonumber\\
= \sum b_{i}\int_{\C}
{\psi}_{\bar z}\frac{a(a-1)}{z(z-1)} \left( \frac{1}{a-c_i} \right) \left(\frac{1}{z-a}-\frac{1}{z-c_i} \right) dz \wedge d{\bar z} 
- \sum \int_{\C}{\psi}_{\bar z}p_{i}(z){\gamma_{a}}(z)dz \wedge d{\bar z}
+\nonumber\\
+ \int_{\C}{\psi}_{\bar z}h(z){\gamma_{a}}(z)dz \wedge d{\bar z}=
\nonumber\\
=\sum \frac{b_i}{a-c_i} \left( \int_{C}{\psi}_{\bar z}{\gamma}_{a}(z)-
\frac{a(a-1)}{c_{i}(c_{i}-1)}\int_{\C}{\psi}_{\bar z}{\gamma}_{c_i}(z) \right) dz \wedge d{\bar z}-\nonumber\\
-\sum \int_{\C}{\psi}_{\bar z}p_{i}(z){\gamma_{a}}(z)dz \wedge d{\bar z} + \int_{\C}{\psi}_{\bar z}h(z){\gamma_{a}}(z)dz \wedge d{\bar z}. 
\end{gather}

On the other hand making some calculations and applying Green's formula we 
have the following equalities.

\begin{displaymath}
\int_{\C}{\psi}_{\bar z}\frac{a(a-1)}{z(z-1)(z-a)}dz \wedge d{\bar z}= 
 (a-1) \left( \int_{\partial supp \varphi} \frac{\psi_{\bar{z}}}{z}  -  a\int_{\partial supp \varphi }\frac{\psi_{\bar{z}}}{z-1} + \int_{\partial supp \varphi} \frac{\psi_{\bar{z}}}{z-a} \right) dz =   
\end{displaymath}

\begin{equation}
= (a - 1)\psi(0) - a(\psi(1)) + \psi(a).
\end{equation}

Since  ${\psi}(0)={\varphi}f(0)={\varphi}(0)=0$, also 
${\psi}(1)={\varphi}f(1)=0$. Hence
${\varphi}(f(a))={\psi}(a) = (3)$.

Applying again Green's formula  we have:

\begin{displaymath}
(3)  = {\varphi}(f(a))={\varphi}(f(a))+(f(a)-1){\varphi}(0)-f(a){\varphi}(1)= 
\int_{\C}\frac{{\varphi_{\bar z}}}{z-f(a)}dz \wedge d{\bar z} + 
\end{displaymath}

\begin{displaymath}
\int_{\C}\frac{(f(a)-1){\varphi_{\bar z}}}{z}dz \wedge d{\bar z}-
\int_{\C}\frac{f(a){\varphi_{\bar z}}}{z-1}dz \wedge d{\bar z}= \int_{\C}{\varphi_{\bar z}}{\gamma_{f(a)}}(z)dz \wedge d{\bar z}
\end{displaymath}
this proves 

\begin{equation}
\int_{\C}{\varphi}_{\bar z}(f){\gamma_{a}}(z)dz \wedge d{\bar z}=\int_{\C}{\varphi_{\bar z}}{\gamma_{f(a)}}(z)dz \wedge d{\bar z}. 
\end{equation}

Applying (4) on (2) we obtain

\begin{gather}
(2) = \sum \frac{b_i}{a-c_i} \left( \int_{\C}{\varphi}_{\bar z}{\gamma}_{f(a)}(z)-\frac{a(a-1)}{c_{i}(c_{i} - 1)}\int_{\C}{\varphi}_{\bar z}{\gamma}_{f(c_i)}(z) \right) dz \wedge d{\bar z} -\nonumber\\
- \sum \int_{\C}{\varphi}_{\bar z}(f)p_{i}(z){\gamma_{f(a)}}(z)dz \wedge d{\bar z}+ \int_{\C}{\varphi}_{\bar z}(f)h(z){\gamma_{f(a)}}(z)dz \wedge d{\bar z} =
\end{gather}

\begin{displaymath}
\int_{\C}{\varphi}_{\bar z}{\gamma}_{f(a)}(z) \left( \sum \frac{b_i}{a-c_i}-p_{i}(a)+h(a) \right) dz \wedge d{\bar z} +
\sum \frac{a (a-1)b_{i}}{c_{i}(c_{i}-1)(c_{i}-a)}\int_{\C}{\varphi}_{\bar z}{\gamma}_{f(c_{i})}(z)dz \wedge d{\bar z} =
\end{displaymath}

\begin{displaymath}
\int_{\C}{\varphi}_{\bar z}{\gamma}_{f(a)}(z)\left( \frac{1}{f'(a)} \right) dz \wedge d{\bar z} + \sum \frac{a(a-1)b_{i}}{c_{i}(c_{i}-1)(c_{i}-a)}\int_{\C}{\varphi}_{\bar z}{\gamma}_{f(c_{i})}(z)dz \wedge d{\bar z}. 
\end{displaymath}

Since ${\gamma}_{a}(c_{i}) = \frac{a(a-1)}{c_{i}(c_{i}-1)(c_{i}-a)}$, we have

\begin{displaymath}
\int_{\C}{\varphi}_{\bar z}f^{*}{\gamma}_{a}dz \wedge d{\bar z} =\int_{\C}{\varphi}_{\bar z} \left[ \frac{{\gamma}_{f(a)}(z)}{f'(a)}+\sum {b_i}{\gamma}_{a}(c_{i}){\gamma}_{f(c_{i})}(z) \right] dz \wedge d{\bar z}.
\end{displaymath}

Hence

\begin{displaymath}
\int_{\C}{\varphi}_{\bar z}(f^{*}({\gamma}_{a}(z))- \left[ \frac{{\gamma}_{f(a)}(z)}{f'(z)}+\sum {b_i}{\gamma}_{a}(c_{i}){\gamma}_{f(c_{i})}(z) \right] dz \wedge d{\bar z}=0. 
\end{displaymath}

This is true for each ${\varphi}$, so the function inside the integral 
is by Weyl's lemma an holomorphic function on $\C \setminus \{0,1\}$ which
is integrable. In our case, this implies that

\begin{displaymath}
f^{*}{\gamma_{a}}(z)= \frac{{\gamma_{f(a)}(z)}}{f'(a)}+\sum_{i =1}{b_i}{\gamma_{a}(c_i)}{\gamma_{f(c_i)}}(z).
\end{displaymath}

\section{Formal Relations of Ruelle Poincare Series} 

In this section we want to study properties of series of the form

\begin{displaymath}
{\sum}_{n=0}^{\infty}x^{n}f^{*n}({\gamma}_{a}(z))
\end{displaymath}

where $f^{*n}$ denotes the n-th iteration of the Ruelle operator. Observe 
from  Section 2 that

\begin{displaymath}
f^{*0}({\gamma}_{a}(z))={\gamma}_{a}(z) 
\end{displaymath}

\begin{displaymath}
f^{*}({\gamma}_{a}(z))=\frac{1}{f'(a)}{\gamma}_{f(a)}(z)+{\sum} b_{i}{\gamma}_{a}(c_{i}){\gamma}_{f(c_{i})}(z)
\end{displaymath}

\begin{displaymath}
f^{*2}({\gamma}_{a}(z))= \frac{1}{f'(a)}(\frac{{\gamma}_{f^{2}(a)}(z)}{f'(f(a))}+ {\sum} b_{i}{\gamma}_{f(a)}(c_{i}){\gamma}_{f(c_{i})}(z)+
{\sum}b_{i}{\gamma}_{a}(c_{i})f^{*}{\gamma}_{f(c_{i})}(z)=
\end{displaymath}

\begin{displaymath}
\frac{\gamma_{f^{2}(a)}(z)}{(f^2)'(a)} + {\sum}b_{i} \left( \frac{1}{f'(a)}
{\gamma}_{f(a)}(c_{i} {\gamma}_{f(c_{i})}(z)+{\gamma}_{a}(c_{i})f^{*}
{\gamma}_{f(c_{i})} \right) 
\end{displaymath}

\begin{displaymath}
f^{*3}({\gamma}_{a}(z))= \frac{1}{(f^3)'(a)}{\gamma}_{f^{3}(a)}(z)+ 
\end{displaymath}

\begin{displaymath}
{\sum} b_{i} \left( \frac{{\gamma}_{f^{2}(a)(c_{i})}{\gamma}_{f(c_{i})}(z)}{(f^2)'(a)} + \frac{{\gamma}_{f(a)}(c_{i})}{f'(a)}f*({\gamma}_{f(c_{i})}(z))+
{\gamma}_{a}(c_{i})f^{*2}({\gamma}_{f(c_{i})}(z) \right)
\end{displaymath}

in general we have

\begin{displaymath}
f^{*n}({\gamma}_{a}(z))=\frac{1}{(f^{n})'(a)}{\gamma}_{f^{n}(a)}(z)+ 
{\sum}_{i} b_{i}c_{n-1}^{i}
\end{displaymath}

\noindent for some coefficients $c_{j}^{i}$, determined by
the {\it Cauchy's product} of two series $A={\sum}a_{i}$ and $B={\sum}b_{i}$
where $C=A{\bigotimes}B={\sum}c_{n}$ and $c_{n}={\sum}_{i=0}^{n}a_{i}b_{n-i}$.

Now define $S(a,z)={\sum}_{n=0}^{\infty} f^{*n}({\gamma}_{a}(z))$, $A(a,z)={\sum}_{n=0}^{\infty} \frac{1}{(f^n)'(a)} \gamma_{f^n(a)}(z)$ and 
$S(f(c_{i}),z){\bigotimes}A(a,z)=C^{i}= {\sum}_{j}c_{j}^{i}$. Thus  we have

\begin{displaymath}
S(a,z)=A(a,z)+ {\sum}b_{i}{\sum}_{n=0}^{\infty}c_{n-1}^{i}=
\end{displaymath}

\begin{displaymath}
A(a,z)+ {\sum}b_{i}[S(f(c_{i}), z){\bigotimes}A(a,c_{i})].
\end{displaymath}

Define $S(x,a,z)={\sum}x^{n}f^{*n}({\gamma}_{a}(z))$ and
$A(x,a,z)={\sum} \frac{x^{n}}{(f^n)'(a)}{\gamma}_{f^{n}(a)}(z)$, since

\begin{displaymath}
x^{n}f^{*n}({\gamma}_{a}(z))=\frac{x^{n}}{(f^n)'(a)}{\gamma}_{f^{n}(a)}(z)+ x{\sum} b_{i}c_{n-1}^{i}x^{n-1}
\end{displaymath}

\noindent then $S(x,a,z)=A(x,a,z)+ x{\sum}b_{i}{\sum}_{n=1}^{\infty}c_{n-1}^{i}x^{n-1}=A(x,a,z)+ {\sum}b_{i}[S(x,f(c_{i}),z){\bigotimes}A(x,a,c_{i})]$
by the Cauchy's  lemma on power series formula it can be written as
$S(x,a,z)=A(x,a,z)+{\sum}b_{i}[S(x,f(c_{i}),z)A(x,a,c_{i})]$  for all  $x$ in 
the disc  of convergence of the series.

\begin{lemma}
For all  $|x|<1$, $S(x,a,z) \subset L_{1}(\C)$.
\end{lemma}

{\it Proof.}
\begin{displaymath}
\int|S(x,a,z)| \leq {\sum}|x^{n}| \int |f^{*n}({\gamma}_{a}(z))|=
{\sum}|x^{n}|||f^{*n}({\gamma}_{a}(z))|| \leq ||{\gamma}_{a}(z)||{\sum}|x^{n}|=\frac
{||{\gamma}_{a}(z)||}{1-|x|}< {\infty}.
\end{displaymath}

\begin{lemma}
If  $a$ is  summable with $a\in \C$, then  $A(x,a,z) \in L_{1}(\C)$ for all $|x|<1$.
\end{lemma}

{\it Proof.} 

\begin{displaymath}
\int|A(x,a,z)|\leq{\sum}
\frac{|x^{n}|}{|(f^n)'(a)|}||{\gamma}_{f^{n}(a)}(z)||, 
\end{displaymath}

\noindent now  by the properties of potential function we have

\begin{displaymath}
\int |{\gamma}_{t}(z)| \leq K|t| |ln|t||
\end{displaymath}

\noindent hence 

\begin{displaymath}
\int|A(x,a,z)| \leq K{\sum}\frac{|x^{n}|}{|(f^n)'(a)|}|f^{n}(a)||ln f^{n}(a)|
\leq K max_{y \in \bigcup f^{n}(a)}|y||ln|y||{\sum} \left| \frac{x^n}{(f^n)'(a)} \right| < {\infty}, 
\end{displaymath}

\noindent if the series $\{f^{n}(a)\}$ is bounded. If not apply that the series
${\sum}\frac{f^{n}(a) |ln |f^{n}(a)||}{(f^n)'(a)}$ absolutely converges.
Then $A(x,a,z) \in L_{1}(C)$.

\bigskip

\begin{corollary}
Under conditions of Lemma 3 we have  $lim_{x \rightarrow 1}||A(x,a,z)- A(1,a,z)||=0$ in $L_{1}(\C)$.
\end{corollary}

{\it Proof.} 
Observe that we can choose $N$ such that $2{\sum}_{i \geq N}|\frac{1}{(f^n)'(a)}| \leq {\epsilon}/2$, let ${\delta}$ such that
$|1-x|<{\delta}$. We have that $|1-x^{N}|{\sum}_{i<N}|\frac{1}{(f^n)'(a)}|
\leq {\epsilon}/2$. Hence

\begin{displaymath}
{\sum}\left| \frac {x^{n}-1}{(f^n)'(a)} \right| \leq {\sum}_{n<N} \left| \frac{x^{n}-1}{(f^n)'(a)} \right| + {\sum}_{n \geq N} \left| \frac{x^{n}-1}{(f^n)'(a)} \right| \leq {\epsilon}/2+2{\sum}_{n \geq N}\left| \frac{1}{(f^n)'(a)}\right| \leq {\epsilon}
\end{displaymath}

so $lim_{x \rightarrow 1}{\sum}|\frac{x^{n}-1}{(f^n)'(a)}|=0$, but

\begin{displaymath}
lim_{x \rightarrow 1}||A(x,a,z)- A(1,a,z)|| \leq K max_{y \in \bigcup f^{n}(a)}(|y||ln|y||){\sum}\frac{x^n-1}{(f^n)'(a)},
\end{displaymath}

\noindent if sequence  $\{f^{n}(a)\}$ is bounded, otherwise use the absolute convergence
of ${\sum}\frac{f^{n}(a) \ln ((f^n)'(a)}{(f^n)'(a)}$, hence the corollary is proved.

\begin{lemma}
If ${\sum}\frac{1}{(f^n)'(a)}$ converges absolutely, then $\lim_{x \rightarrow 1}|A(x,a,c_{i})- A(a,c_{i})|=0$.
\end{lemma}

{\it Proof.}  Consider  $X_{a}={\bar {\bigcup}_{n>0}f^{n}(a)}$. If $c_{i}$ is not in $X_{a}(f)$, then 
$A(a,c_{i})={\sum}\frac{1}{(f^n)'(a)}{\gamma}_{f^{n}(a)}(c_{i})$ has no poles. 
Since $|{\gamma}_{f^{n}(a)}(c_{i})|<1/d^3$, for some constant $d$, it is 
bounded in $Y_a$.

If $c_{i}$ is in $X_a(f)$, let $D_{\epsilon}=|z-c_{i}|<{\epsilon}$ so for $z\in D_{\epsilon}$, we can use the equality 
$f'(z)= (z-c_{i})f''(c_{i})+O(|z-c_{i}|^2)$ and obtain

$$ \frac{1}{\vert f^{n_i} - c_i\vert} \leq \frac{\vert f''(c_i))\vert + O(\vert f^{n_i} - c_i\vert)}{\vert f'(f^{n_i}(a))\vert} \leq K\frac{1}{\vert f'(f^{n_i}(a))\vert}, 
$$
hence
$$
{\biggl |}\gamma_{f^{n_i}(a)}(c_i){\biggr |} \leq K\frac{1}{\vert f'(f^{n_i}(a))\vert}\frac{\vert f^{n_i}(a)\vert\vert f^{n_i}(a) - 1\vert}{\vert c_i(c_i - 1)} \leq K_1 \frac{1}{\vert f'(f^{n_i}(a))\vert}, 
$$
where $K$ and $K_1 $ are constant depending only on $\epsilon $ and the points 
$c_i$. As result for all $ \vert x\vert \leq 1 $ we have

$$ {\biggl |}\sum_i\frac{x^{n_i}}{(f^{n_i})'(a)}\gamma_{f^{n_i}(a)}(c_i){\biggr |} \leq K_1\sum_i\frac{\vert x\vert^{n_i}}{\vert (f^{n_i + 1})' (a) \vert} < \infty.
$$

This proves the Lemma. So we have that the following 
equality holds

\begin{displaymath}
S(x,a,z)=A(x,a,z)+x{\sum}b_{i}S(x,a,d_{i})A(x,a,c_{i}).
\end{displaymath}


\section{Ruelle Operator and Line Fields}

Let $f$ be a transcendental entire map, we say that $f$ admits an {\it  invariant line field} if there is a measurable Beltrami differential ${\mu}$ on the 
complex plane $\C$ such that $B_{f}{\mu}={\mu}$ a.e. $|{\mu}|=1$ on a set of
 positive measure and ${\mu}$ vanishes else were. If ${\mu}=0$
outside the Julia set $J(f)$, we say that ${\mu}$ is 
{\it carried on the Julia set}. See \cite{McM} for results of holomorphic line 
fields.

\bigskip

In Section 2 we consider the set $Fix (B_{f})=\{{\mu} \in L_{\infty}(C):B_{f}({\mu})={\mu}\}$, with $B_{f}$ 
being the Beltrami operator. Consider now the following integrals

\begin{displaymath}
\int_{\C}{\mu}S(x,a,z)dz \wedge d{\bar z}= \int{\mu}A(x,a,z)dz \wedge d{\bar z}+ x{\sum}b_{i}A(x,a,c_{i})
\int{\mu}S(x,f(c_{i}),z)dz \wedge d{\bar z}.
\end{displaymath}

The  above equation is  equal to the following expression, by 
the properties of the potential $F_{\mu}$

\begin{displaymath}
{\sum}x^{n}\int{\mu}(f^{*})^n{\gamma}_{a}(z)=
{\sum}\frac{x^{n}}{(f^n)'(a)}F_{\mu}(f^{n}(a))+x{\sum}b_{i}A(x,a,c_{i}) \int \mu \sum x^n (f^*)^n(\gamma_{f(c_i)}(z)) =
\end{displaymath}

\begin{displaymath}
{\sum}\frac{x^{n}}{(f^n)'(a)}F_{\mu}(f^{n}(a))+ \frac{x}{1 - x} {\sum}b_{i}A(x,a,c_{i})F_{\mu}(f^n(c_{i})).
\end{displaymath}

By invariance of the Ruelle operator we have

\begin{displaymath}
\sum x^n \int \mu (f^*)^n {\gamma}_{a}(z) =\frac{\int{\mu}{\gamma}_{a}(z)}{1-x}=\frac{F_{\mu}(a)}{1-x}=
\end{displaymath}

\begin{displaymath}
{\sum}\frac{1}{(f^n)'(a)}x^{n}F_{\mu}(f^{n}(a))+\frac{x}{1-x}{\sum}b_{i}A(x,a,c_{i})F_{\mu}(f^n(c_{i})).
\end{displaymath}

Hence 

\begin{equation}
F_{\mu}(a) =(1-x) \sum \frac{1}{(f^n)'(a)} x^{n} F_{\mu}(f^{n}(a)) + x \sum b_{i}A(x,a,c_{i})
F_{\mu} (f^n(c_{i})).
\end{equation}

By Corollary 2 and Lemma 4 we can pass to the limit  $x \to 1$  in (5), as a result we have
 
\begin{displaymath}
F_{\mu}(a) = \sum b_i A(1, a, c_i) F_{\mu}(f^n(c_{i})).
\end{displaymath}

By hypothesis $d_1$ is summable, so for $f \in S_{q}$ 

\begin{eqnarray}
F_{\mu}(d_{1}) \left( 1 - \sum_{\genfrac{}{}{0pt}{}{c_k}{f(c_k) = d_1}} b_k A(d_1 ,c_k) \right) = \sum_{i=2}^{q} F_{\mu}(d_{i}) \left(  \sum_{\genfrac{}{}{0pt}{}{c_l}{f(c_l) = d_i}} b_l A(d_{1}, c_l) \right). \nonumber
\end{eqnarray}

If we denote $\Psi_i = \sum_{c_k} b_k  A(d_1 ,c_k)$ for $f(c_k) = d_i$, 
then we can rewrite the above equation as:

\begin{equation}
F_{\mu}(d_{1})(1 - \Psi_1) = \sum_{i = 2}^{q} F_{\mu}(d_{i}) \Psi_i.
\end{equation}

\begin{definition}
We say that  (6) is a trivial relation if and only if $\Psi_i = 0$, 
$i = 2, 3 \dots q$  and  $\Psi_1 = 1$. 
\end{definition}
 
For transcendental entire maps  there are, in general, many critical 
points $c_{k}$ which are mapped to the critical value $d_{1}$, even if the 
function is structurally stable.

\begin{proposition} If (6) is   a non  trivial relation, then $f$ is 
unstable.
\end{proposition}

{\it Proof.} By Hypothesis the set of $\{ f(c_i) \}$ is finite. By 
equation (5)  the Bers operator $\beta$ induces an isomorphism

$$
\beta^* : HD(f) \times  J_f \to \C^q.
$$

with  coordinates 
$\beta^*(\mu) = \{ F_{\mu}(d_{1}) \dots  F_{\mu}(d_{q}) \}$.

If (6) is  a non trivial relation, then the relation gives a non trivial 
equation on the image of $\beta^*$, where the  image of  $\beta^*$ is a 
subset of the set of solutions of this equation. Then  
$dim( HD(f) \times  J_f) = dim(\text{image of } \beta^*) < q$. Thus 
the proposition is proved.


\section{Fixed Point Theory}

In this section we want to prove Theorems A and B  which were stated in the 
introduction. In order to prove the theorems we will give a series of results.

\begin{proposition}
If (6) is a trivial relation, then $f^{*}A(d_{1},z)=A(d_{1},z)$.
\end{proposition}
\bigskip

{\it Proof.}
Let us remember that $A(d_{1},z)={\sum}\frac{1}{(f^n)'(d_{1})}{\gamma}_{f^{n}(d_{1})}(z)$
and

\begin{displaymath}
f^{*}({\gamma}_{a}(z))=\frac{1}{f'(a)}{\gamma}_{f(a)}(z)+ {\sum} b_{i}{\gamma}_{a}(c_{i}){\gamma}_{f(c_{i})}(z), 
\end{displaymath}

\noindent and so 

\begin{displaymath}
f^{*}(A(d_{1},z))= A(d_{1},z)-{\gamma}_{d_{1}}(z)+{\sum}_{i}{\gamma}_{f(c_{i})}(z) \left( \sum_{\genfrac{}{}{0pt}{}{c_l}{f(c_l) = d_i}} b_l A(d_1, c_l) \right) = A(d_1, z)
\end{displaymath}

\medskip

Denote by $Z=\overline{\bigcup_{i}f^{i}(d_{1})}$ and $Y=\C \setminus Z$

\medskip
 
\begin{proposition} 
If $\varphi = A(d_1, z) \neq 0$ on $Y$, then (6) is a non trivial relation.
\end{proposition}

 Before we prove the above proposition we will prove a series of results 
which will help us to prove the proposition.

\bigskip

Consider the modulus of the  Ruelle operator: $|f^{*}|{\alpha}={\sum}{\alpha}{(\xi}_{i})|{\xi}_{i}'|^2$ where ${\xi}_{i}$ are the inverse branches of $z$ under the map $f$.

\begin{lemma}
$|f^{*}||{\varphi}|=|{\varphi}|$.
\end{lemma}

By hypothesis we have ${\varphi}=f^{*}{\varphi}={\sum}{\varphi}({\xi}_{i}){(\xi}_{i}')^2$ For fixed $i$, denote by ${\alpha}_{i}(z)={\varphi}({\xi}_{i}){(\xi}_{i}')^2$ and ${\beta}_{i}= {\varphi}-{\alpha}_{i}$. We have the following 
claim:
\vspace{.5cm}

\noindent {\bf Claim.} {\it $|{\alpha}_{i}+{\beta}_{i}|=|{\alpha}_{i}|+|{\beta}_{i}|$,  for almost every point.}
\vspace{.5cm}

{\it Proof.}

\begin{displaymath}
||{\varphi}||=||f^{*}{\varphi}||=\int|{\alpha}_{i}+{\beta}_{i}| \leq \int|{\alpha}_{i}|+ \int||{\beta}_{i}|\leq ||{\varphi}||, 
\end{displaymath}

\noindent which  implies that $\int|{\alpha}_{i}+{\beta}_{i}| = \int|{\alpha}_{i}|+\int|{\beta}_{i}|$. Now let $A=\{z:|{\alpha}_{i}(z)+{\beta}_{i}(z)| < |{\alpha}_{i}(z)|+|{\beta}_{i}(z)| \}$ with $m(A) \geq 0$, where  $m$ is  the Lebesgue measure. Then $\int_{(\C \setminus A)\bigcup A}|{\alpha}_{i}+{\beta}_{i}|  \int_{A}|{\alpha}_{i}+{\beta}_{i}|+
\int_{\C-A}|{\alpha}_{i}+{\beta}_{i}|< \int_{A}|{\alpha}_{i}(z)|+|{\beta}_{i}(z)|+\int_{\C \setminus A}|{\alpha}_{i}(z)|+|{\beta}_{i}(z)|$, which is a contradiction, thus the claim is proved.

\bigskip

Now by induction on the claim, we have that ${\sum}|{\alpha}_{i}|=|{\sum}{\alpha}_{i}|$ and so $f^{*}|{\varphi}|=|{\varphi}|$. This proves Lemma 5.
\bigskip

\noindent{\bf Remark 2.} The measure $\sigma(a) = \int \int_A |\phi(z)|$ is a 
non negative invariant absolutely continue probability measure, where 
$A \subset \chat$ is a measurable set.
\bigskip

\begin{definition}
A measurable set $A \in \chat$ is called back wandering if and only if $m(f^{-n}(A) \cap f^{-k}(A)) = 0$, for $ k \neq n$. 
\end{definition}

\begin{corollary}
If $\varphi \neq 0$ on $Y$, then (i) $J(f) = \chat$, (ii) $m(Z) = 0$ and  
(iii) $\frac{\bar{\varphi}}{\vert \varphi \vert}$  defines an invariant Beltrami differential. 
\end{corollary}

{\it Proof.} (i) Every non periodic point of the Fatou set has a back 
wandering 
neighborhood.  By Remark 2  we have that $\varphi  = 0$. Thus 
$J(f) = \chat$ and  $\varphi \neq  0$ on every component of $Y$.

\bigskip

(ii) If $m(Z) > 0$,  then $m(f^{-1}(Z)) >0$ so $m(f^{-1}(Z) - Z) >0 $ since
$f^{-1}(Z) \neq Z$, $Z \neq \C$, denote by $Z_1 =  f^{-1}(Z) - Z$. Then   
$Z_1$ is  back wandering thus $\varphi = 0$. Therefore, $m(Z) = 0$. 

\bigskip

(iii) By using notations and the proof of Lemma 5 we have 
  $f_i(x)= \frac{{\alpha}_{i}}{{\beta}_{i}}= \frac{{\varphi}}{{\alpha}_{i}}-1$
so ${\varphi}=(1+f_{i}(x)){\alpha}_{i}=(1+f_{i}(x)(\varphi(\xi_i(x))(\xi_i')^{2}(x)$, 
with $x \in supp({\alpha}_{i})$. Consider $t={\xi}_{i}(x)$, with $t \in {\xi}_{i}(supp ({\alpha}_{i}))$. 
Then ${\varphi}(f(t))(f')^{2}(t)=(1+f_{i}(f(t))){\varphi}(t)$. Hence, 

\begin{displaymath}
\frac{{\bar {\varphi}(x)}}{|{\varphi}(x)|}=\frac{(1+f_{i}(x)){\bar{\varphi}({\xi}_{i}(x))}{\bar{(\xi_i')}^{2}(x)}}{(1+f_{i}(x))|{\varphi}({\xi}_{i}(x)||(\xi_i')^{2}(x)|}, 
\end{displaymath}

\noindent and so
 
\begin{displaymath}
{\mu}=\frac{{\bar {\varphi}}}{|{\varphi}|}=\frac{{\bar {\varphi}({\xi}_{i})}{\bar {\xi}_{i}'}}{|{\varphi}({\xi}_{i})|{\xi}_{i}'}
\end{displaymath}

as result $\mu = \mu (\xi_i)\frac{\bar{\xi_i'}}{\xi_i'}$ is an invariant line 
field. Thus the corollary is proved.

\bigskip 

\begin{lemma} $\frac{\beta_j}{\alpha_j} = k_j \geq 0 \text{ is a non negative function.}$
\end{lemma}

{\it Proof}. We have $ \vert 1 + \frac{\beta_j}{\alpha_j}\vert = 1 + {\biggl |}\frac{\beta_j}{\alpha_j}{\biggr |}, $ then if $ \frac{\beta_j}{\alpha_j} = \gamma^j_1 + i\gamma^j_2 $ we have
$$ {\biggl (}1 + (\gamma^j_1){\biggr )}^2 + (\gamma^j_2)^2 = {\biggl (}1 + \sqrt{(\gamma_1^j)^2 +(\gamma_2^j)^2}{\biggr )}^2 = 1 + (\gamma_1^j)^2 +(\gamma_2^j)^2 + 2\sqrt{(\gamma_1^j)^2 +(\gamma_2^j)^2}.
$$
Hence $ \gamma_2^j = 0  $ and $ \frac{\alpha_j}{\beta_j} = \gamma_1^j $ is a real-valued function but $\frac{\alpha_j}{\beta_j}$ is meromorphic function. So $ \gamma_1^j = k_j $ is constant on every connected component of $ Y $ and the condition $ \vert 1 + k_j\vert = 1 +\vert k_j\vert $ shows $ k_j \geq 0.$

\bigskip

\noindent {\bf Proof of Proposition 6.4.}

{\it Proof.} Let us show first that all postcritical values are in Z. 
Assume that there is some $d_{i} \in Y$, then by the  Lemma 6,   ${\varphi}(z)=(1+k_{j}){\varphi}({\xi}_{j}(z)){\xi'}_{j}^{2}(z)$.
Assume that the branch ${\xi}_{j}(z)$ is such that tends to $c_{i}$ when 
$z$ tends to $d_{i}$. Then ${\xi'}_{j}^{2}$ tends to $\infty$ and so ${\varphi}(c_{i})=0$.
Also for every $k_{j}$ we have that ${\varphi}(c_{i})=(1+k_{j}){\varphi}({\xi}_{j}(c_{i})){\xi'}_{j}^{2}(c_{i})$, so ${\varphi}({\xi}_{j}(c_{i}))=0$. This implies that if $c$ is a preimage of a critical point,
then ${\varphi}(c)=0$, since $J=\C$ then ${\varphi}=0$ in $Y$, which is a contradiction.

Let us show now that $Z={\bigcup}f^{i}(d_{1})$. 
We will use a McMullen argument like in \cite{McM}. By Lemma 5 and Corollary 3,
$\mu=\frac{\bar{\varphi}}{|{\varphi}|}$ is an invariant line field. That implies that ${\varphi}$ is dual to ${\mu}$ and
it is defined up to a constant. We will construct a meromorphic function 
$\psi$, dual to $\mu$ and such that $\psi$ has finite number of poles on each disc $D_{R}$ of radius $R$ centered at $0$. 

For that suppose that for $z \in \C$ there exists a branch
$g$ of a suitable $f^{n}$, such that $g(U_{z}) \in Y$, where $U_{z}$ is a neighborhood of $z$.
Then define $\psi(\xi)=\varphi(g(\xi))(g')^{2}(\xi)$, for all $\xi \in U_{z}$. Note that $\psi(\xi)$ is 
dual to $\mu$ and has no poles in $U_{z}$. If there is no such branch $g$, then $\xi$ is in the postcritical set, and there is a branch covering $F$ from a neighborhood of $\xi$ to $U_{z}$,
then define $\psi(\xi)=F^{*}(\varphi)$, with $F^{*}$ the Ruelle operator of $F$. The map $\psi$ is a meromorphic function dual to $\mu$ in $U_{z}$ and has finite number of poles.

By the discussion above it is possible to construct a meromorphic function dual to $\mu$ in 
any compact disc $D_{R}$. If we make $R$ tends to $\infty$, we have a meromorphic function
$\psi$ defined in $\C$ dual to $\mu$ and with a discrete set of poles. 

Observe now that such $\psi$ is holomorphic in $Y$, then $Z$ is discrete and so $Z={\bigcup}f^{i}(d_{1})$ as we claim. Since every postcritical set is in $Z$, that implies that $f$ is unstable and this proves the proposition. 

\bigskip

The following propositions can be found in \cite{Mak2}. For completeness  we
 prove them.

\begin{proposition} Let $ a_i \in \C, a_i\neq a_j, $ for $ i\neq j $ be points such that $ Z = \overline{\cup_i a_i} \subset {\C} $ is a compact set. Let  $ b_i \neq 0 $ be complex numbers such that the series $ \sum b_i $ is absolutely convergent. Then the function $ l(z) = \sum_i\frac{b_i}{z - a_i} \neq 0 $ identically on $ Y = {\C} \backslash Z $ in the following cases
\begin{enumerate}
\item the set $ Z $ has zero Lebesgue measure
\item if diameters of components of $\C \backslash Z $ uniformly bounded below from zero and
\item If $ O_j $ denote the components of $ Y, $ then $ \cup_i a_i \in \cup_j{\partial O_j}.$
\end{enumerate}
\end{proposition}

{\it Proof}  Assume that $ l(z) = 0 $ on $ Y.$ Let us calculate derivative $ \overline{\partial}l $ in sense of distributions, then $ \omega = \overline{\partial}l = \sum_i b_i\delta_{a_i} $ and by standard arguments 

$$ l(z) = - \int\frac{d\omega(\xi)}{\xi - z}.
$$

Such as $ a_i \neq a_j $ for $ i \neq j, $  then measure $ \omega =0 $ iff all coefficients $ b_i = 0. $ 

Let us check (1). Otherwise in this case we have that the function $ l $ is locally integrable and $ l = 0 $ almost everywhere and hence $ \omega = \overline{\partial}l = 0 $ in sense of distributions and hence $ \omega = 0 $ as a functional on space of all continuous functions on $ Z $ which is a contradiction with the arguments above.

2)   Assume that $ l = 0 $ identically out of $ Z.$ Let $ R(Z) \subset C(Z) $ denote the algebra of all uniform limits of rational functions with poles out of $ Z $ in the $\sup-$topology, here $ C(Z) $ as usually denotes the space of all continuous functions on $ Z $ with the $ \sup-$norm. Then measure $ \omega $ denote a lineal functional on $ R(Z).$ 
The items (2) and (3) are based on the generalized Mergelyan theorem (see \cite{ Gam}) which states {\it If diameters of all components of $ \C \backslash Z $ are bounded uniformly below from 0, then  every continuous function holomorphic on interior of $ Z $ belongs to $ R(Z). $}

Let us show that $ \omega $ annihilates the space $ R(Z). $ Indeed let $ f(z) \in R(Z) $ be a transcendental entire map  and $ \gamma $ enclosing $ Z $ close enough to $ Z $ such that $ f(z) $ does not have poles in interior of $ \gamma.$ Then such that $ l = 0 $ out of $ Z $ we only apply Fubini's theorem
$$ 
\int r(z)d\omega(z) = \int d\omega(z)\frac{1}{2\pi i}\int_{\gamma}\frac{r(\xi)d\xi}{\xi - z} = \frac{1}{2\pi i}\int_{\gamma}r(\xi)d\xi\int\frac{d\omega(z)}{\xi - z} = \frac{1}{2\pi i}\int_{\gamma}r(\xi)l(\xi)d\xi = 0.
$$
Then by generalized Mergelyan theorem we have $ R(Z) = C(Z) $ and $ \omega = 0. $ Contradiction.

\bigskip

Now let us check (3). We {\bf claim} {\it that $ l = 0 $ almost everywhere on} $ \cup_i\partial O_i. $

\bigskip 

{\it Proof of the claim.}  Let $ E \subset \cup_i\partial O_i $ be any measurable subset with positive Lebesgue measure. Then the function $ F_E(z) = \iint_E\frac{dm(\xi)}{\xi - z} $ is continuous on $\C  \backslash \cup_i O_i $ and is holomorphic onto interior of $\C \backslash  \cup_i O_i. $ Again by generalized Mergelyan theorem  $ F_E(z) $ can be approximated on $\C  \backslash \cup_i O_i $ by functions from $ R(\C \backslash  \cup_i O_i) $ and hence by arguments above and by assumption we have
$ \int F_E(z)d\omega(z) = 0. $ But again application of Fubini's theorem gives

$$ 0 = \int F_E(z)d\omega(z) = \int d\omega(z)\iint_E\frac{dm(\xi)}{\xi - z} = \iint_E dm(\xi)\int d\omega(z)\frac{1}{\xi - z} = \iint l(\xi)d(m(\xi)).
$$
Hence for any measurable $ E \subset \cup_i\partial O_i $ we have $ \iint_E l(z) = 0. $ The claim is proved.
Now for any component $ O \in Y $ and  any measurable $ E \subset \partial O $ we have $ \iint_E l(z) = 0. $ By assumption $ l = 0 $ almost everywhere on 
$\C$. Contradiction thus the proposition is proved.

\begin{proposition} If $f \in X$, then $ A(d_1 , z) = \varphi(z) \neq 0 $ identically on $ Y$ in the following cases 
\begin{enumerate}
\item if $ f^{-1}(d_{1}) \notin X_{c_1}, $
\item if diameters of components of $ Y $ are uniformly bounded below from 0, 
\item  If $ m\left (X_{c_1}\right ) = 0, $ where $ m $ is the Lebesgue measure on $ {\C}, $
\item if $ X_{c_1} \subset \cup_i \partial D_i, $ where $ D_i $ are components of Fatou set. 
\end{enumerate}
\end{proposition}

{\it Proof} Let us prove (1). If $ f $ is structurally stable then relation 
(3) is trivial. 

Assume now that the set $ X_{c_1} $ is bounded. Then by Proposition 6.5  
we have that $ \varphi(z) = \frac{C_1}{z}+\frac{C_2}{z-1}+ \sum\frac{1}{(f^i)'(d_1)(z - f^i(d_1))} = l(z) \neq 0$.  Other cases follows directly from Proposition 6.5 also.

\par Now let $ X_{c_1} $ be unbounded. Let $ y \in {\C} $ be a point such that the point $ 1 - y \in Y, $ then the map $ g(z) = \frac{y z}{z + y - 1} $ maps $ X_{c_1} $ into $ \C. $ Let us consider the function $ G(z) = \frac{1}{z}\sum_i\frac{(f^i(d_1) - 1)}{(f^i)'(d_1)} - \frac{1}{z - 1}\sum_i\frac{f^i(d_1)}{(f^i)'(d_1)} +  \sum\frac{1}{(f^i)'(d_1)(z - g(f^i(d_1))}, $ then by proposition 6.5 $ G(z) \neq 0 $ identically on $ g(Y). $

 Now we {\bf Claim} that {\it Under condition of theorem A }
$$
 G(g(z))g'(z) = \phi(z).
$$ 
\par {\it Proof of claim.} Let us define $ C_1 = \sum_i\frac{(f^i(d_1) - 1)}{(f^i)'(d_1)} $ and $ C_2 = \sum_i\frac{f^i(d_1)}{(f^i)'(d_1)} $ then we have
$$
\frac{C_1}{g(z)} = \frac{C_1(z + y - 1)}{yz} \text{ and } \frac{C_2}{g(z) - 1} = \frac{C_2(z + y - 1)}{(y - 1)(z - 1)}
$$
and for any $ n $

\begin{displaymath}
\frac{1}{g(z) - g(f^n(d_1))} = \frac{(z + y - 1)(f^n(d_1) + y - 1)}{y(y - 1)(z - f^n(d_1))} = \frac{1}{y(y - 1)}\left(\frac{(z + y - 1)^2}{z - f^n(d_1)} + 1 - y - z\right),
\end{displaymath}
then 

\begin{multline}
G(g(z)) = \frac{C_1(z + y - 1)}{yz} -  \frac{C_2(z + y - 1)}{(y - 1)(z - 1)} + \sum \frac{1}{(f^i)'(d_1)(g(z) - g(f^i(d_1))} =\nonumber\\ 
= \frac{1}{y(y - 1)}\biggl((1 - y - z)\sum \frac{1}{(f^i)'(d_1)} + (z + y - 1)^2\sum\frac{1}{(f^i)'(d_1)(z - f^i(d_1)} + \nonumber\\
+ \frac{C_1(z + y - 1)}{yz} -  \frac{C_2(z + y - 1)}{(y - 1)(z - 1)}\biggr) = \nonumber\ast
\end{multline}
and 
\begin{gather}
\ast = \frac{1}{g'(z)}\left( \phi(z) - \frac{\sum_i\frac{f^i(d_1) - 1}{(f^i)'(d_1)}}{z} + \frac{\sum_i\frac{f^i(d_1)}{(f^i)'(d_1)}}{z - 1} + \frac{\sum\frac{1}{(f^i)'(d_1)}}{1 - y - z} + \frac{C_1(y - 1)}{z(z + y - 1)} - \frac{C_2y}{(z - 1)(z + y - 1)}\right)= \nonumber\\
= \frac{\phi(z)}{g'(z)}.\nonumber
\end{gather}

Hence $ \phi(z) = 0 $ identically on $ Y $ if and only if  $ G(z) = 0 $ identically  on $ g(Y).$ So by proposition 6.5 we complete this proposition.

\bigskip

\noindent $\underline{\rm {\bf Proof \; \;  of \; \; Theorem \; \; A}}$
\bigskip 

\begin{thmA} Let $f \in X$. If $f$ has a summable 
critical point, then $f$ is not structurally stable map.
\end{thmA}

{\it Proof.} It follows from  Proposition 6.6 that  $\varphi \neq 0$ on $Y$, 
then (6) is non a trivial relation by Proposition 6.4, then applying 
Proposition 5.2 the map $f$ is not stable. Therefore  Theorem A is proved.
\bigskip

\begin{coroA}
Let $f$ transcendental entire map with summable critical point $c \in J(f)$. 
If ${\varphi} \neq 0$ onto $ \C \backslash X_c, $ then $f$ is an unstable map.
\end{coroA}

\noindent $\underline{\rm {\bf Proof \; \;  of \; \; Theorem \; \; B}}$

\bigskip

\begin{thmB} If $f \in W$ is  summable, then there exists no invariant line 
fields on $ J(f).$
\end{thmB}

{\it Proof.} As we observe in equation (6) in Section 4, each summable critical point restricts the image of the 
${\beta}$ operator. The image of the ${\beta}$ operator , belongs
to the common solutions of the equations

\begin{displaymath}
F_{{\mu}}(f(c_{i})(1-b_{i}A(c_{i},f(c_{i}))=\sum_{i \neq j}b_{j}F_{{\mu}}(f(c_{j})A(c_{j},f(c_{j})
\end{displaymath}

\noindent for all $c_{i} \in J(f)$, hence if this system is linearly 
independent,
then the dimension of the image of ${\beta}$ will be $0$ and so 
we will have $J_{f}= \emptyset$. So we have to assume that the above system 
is linearly dependent.

That means in this case, that there are constants
$B_{i}$ such that the function ${\varphi}(z)=\sum B_{i}A(z,f(c_{i}))$ 
is a fixed point of the Ruelle operator $f^{*}$. As in the Lemmas above, the 
measures

$$
\frac{\partial \varphi}{\partial \overline{z}} = \sum_i B_{i} \sum_n \frac{\delta_{f^n (f(c_i))}}{(f^n)'(f(c_i))}=0.
$$

Then $B_i = 0$, this proves Theorem B.

\noindent{\tt Petr M. Makienko}\\
{Permanent addresses:}\\
{ \small Instituto de Matematicas, UNAM}\\ 
{\small Av. de Universidad s/N., Col. Lomas de Chamilpa}\\
{\small Cuernavaca, Morelos, C.P. 62210, M\'exico.}\\
{\small \tt E-mail: makienko@aluxe.matcuer.unam.mx}\\
 {\small and}\\ 
{\small Institute for Applied Mathematics,}\\
{\small   9 Shevchenko str.,}\\
{\small Khabarovsk, Russia}\\
{\small \tt E-mail makienko@iam.khv.ru}

\bigskip

\noindent{\tt Guillermo Sienra}\\ 
{ \small Facultad de Ciencias, UNAM}\\
{ \small Av. Universidad 30, C.U.}\\
{\small M\'exico D.F., C.P. 04510,   M\'exico.}\\
{\small \tt E-mail gsl@hp.fciencias.unam.mx} 

\bigskip

\noindent{\tt Patricia Dom\'{\i}nguez} \\ 
{\small F.C. F\'{\i}sico-Matem\'aticas, B.U.A.P}\\
{\small Av. San Claudio, Col. San Manuel,  C.U.}\\
{\small Puebla Pue., C.P. 72570,  M\'exico}\\
{\small \tt  E-mail pdsoto@fcfm.buap.mx}

\end{document}